\documentclass[12pt]{amsart}
\usepackage{amsmath}
\usepackage{amssymb, amsfonts, amsmath, amsthm, amscd}

\def\deg{\rm{deg}}%
\def\dim{\rm{dim}}%
\def\max{\rm{max}}%

\def\cli{\hbox{\rm Cliff}}

\newtheorem{lem}{Lemma}[section]
\newtheorem{thm}[lem]{Theorem}

\newtheorem{Prop}[lem]{Proposition}
\newtheorem{rmk}[lem]{Remark}
\newtheorem{cor}[lem]{Corollary}
\newtheorem{ex}[lem]{Example}
\newtheorem{question}[lem]{Question}

\begin{document}
\title[Normal generation and Clifford index]
      {Normal generation and Clifford index}
\author[Youngook Choi]{Youngook Choi$^1$}
\address{$^1$Department of Mathematics,
   Seoul National University, Seoul 151-742, Korea}
\email{ychoi@math.kaist.ac.kr}

\author[Seonja Kim]{Seonja Kim$^2$}
\address{$^2$Department of Electronics, Chungwoon University, Chungnam, 350-701, Korea}
\email{sjkim@chungwoon.ac.kr}

\author[Young Rock Kim]{Young Rock Kim$^3$}
\address{$^3$Department of Mathematics, Konkuk University, Seoul, 143-701, Korea}
\email{rocky777@math.snu.ac.kr}

\thanks{This work was supported by Korea Research Foundation Grant funded by Korea Government
(MOEHRD, Basic Reasearch Promotion Fund)(KRF-2005-070-C00005).}

\subjclass[2000]{14H45, 14H10, 14C20}

\keywords{algebraic curve, linear system, line bundle, Clifford
index, normal generation, extremal line bundle.}

\begin{abstract}
  Let $C$ be a smooth curve of genus $g\ge 4$ and Clifford index $c$.
  In this paper, we prove that if $C$ is neither hyperelliptic nor
  bielliptic with $g\ge 2c+5$ and $\mathcal M$ computes the Clifford
  index of $C$, then either $\deg \mathcal M\le \frac{3c}{2}+3$
  or $|\mathcal M|=|g^1_{c+2}+h^1_{c+2}|$ and $g=2c+5$. This
  strengthens the Coppens and Martens' theorem (\cite{CM}, Corollary 3.2.5).
  Furthermore, for the latter case (1) $\mathcal M$ is half-canonical unless $C$ is
  a $\frac{c+2}{2}$-fold covering of an elliptic curve, (2)
  $\mathcal M(F)$ fails to be normally generated with $\cli(\mathcal M(F))=c$, $h^1(\mathcal
  M(F))=2$ for $F\in g^1_{c+2}$.
  Such pairs $(C,\mathcal M)$ can be found on a $K3$-surface whose Picard
  group is generated by a hyperplane section in $\mathbb P^r$. For such a $(C, \mathcal
  M)$ on a K3-surface, $\mathcal M$ is normally generated while $\mathcal M(F)$ fails to be
  normally generated with $\cli(\mathcal M)=\cli(\mathcal M(F))=c$.
\end{abstract}
\maketitle
\pagestyle{plain}

\section{Introduction}
\label{sec:introduction}

Let $C$ be an irreducible projective curve over an algebraically
closed field of characteristic zero. A smooth curve $C$ in $\mathbb
P^r$ is said to be projectively normal if the natural morphisms
$H^0(\mathbb P^r,{\mathcal O}_{\mathbb P^r}(m))\to H^0(C,{\mathcal
O}_C(m))$ are surjective for every nonnegative integer $m$. A line
bundle $\mathcal L$ on a smooth curve $C$ is said to be normally
generated if $\mathcal L$ is very ample and $C$ has a projectively
normal embedding via its associated morphism $\phi_{\mathcal L}:C\to
\mathbb P(H^0(\mathcal L))$.

Green and Lazarsfeld gave a sufficient condition for a line bundle
to be normally generated as follows (\cite{GL}, Theorem 1): If
$\mathcal L$ is a very ample line bundle on $C$ with ${\deg}\mathcal
L\ge 2g+1-2h^1(\mathcal L)-\cli(C)$, then $\mathcal L$ is normally
generated. Note that the condition ${\deg}\mathcal L\ge
2g+1-2h^1(\mathcal L)-\cli(C)$ is equivalent to $\cli(\mathcal
L)<\cli(C)$, whence $h^1(\mathcal L)\le 1$. A very ample line bundle
$\mathcal L$ on a smooth curve $C$ is said to be {\it extremal} if
$\cli(\mathcal L)=\cli(C)$ and $\mathcal L$ fails to be normally
generated.

In this paper, we show that there are extremal line bundles
$\mathcal L$ on smooth curves with $h^1(\mathcal L)= 2$. The
existence of an extremal line bundle $\mathcal L$ with $h^1(\mathcal
L)\le 1$ can be found in \cite{GL}, \cite{Ko};
\begin{thm}[Green-Lazarsfeld, \cite{GL}]
Let $C$ be a smooth curve of genus $g$ and Clifford index $c$. If
$g>{\max} \{\binom{c+3}{2}, 10c+6 \}$, then:

\noindent $(a)$ $C$ always carries an extremal line bundle $\mathcal
L$ with $h^1(\mathcal L)=0$, but never one with $h^1(\mathcal L)\ge
2$;

\noindent $(b)$ $C$ carries an extremal line bundle $\mathcal L$
with $h^1(\mathcal L)=1$ if and only if $c=2f\ge 4$ is even, and $C$
is a two-sheeted branched covering $\pi:C\to C'\subset \mathbb P^2$
of a smooth plane curve $C'$ of degree $f+2$.
\end{thm}

 Thus all the extremal line bundle $\mathcal L$ with $h^1(\mathcal L)\ge
 2$ could be found in case $\cli(C)$ is not so small compared to the genus $g$.
 By the way it is well known that
 $\cli(C)\le [(g-1)/2]$ for any smooth curve and $\cli(C)= [(g-1)/2]$
 if $C$ is a general curve \cite{Me}.
 To find an extremal line bundle $\mathcal L$ with
 $h^1(\mathcal L)\ge 2$, we investigate the property of a line bundle computing
the Clifford index of $C$ and the normal generation of a line bundle
on $C$. We prove the following theorem which supplies a tool
constructing various non-normally generated line bundles on a curve.

\begin{thm}\label{main2}
 Let $\mathcal M$ be a birationally very ample line bundle on a smooth curve $C$,
 $\mathcal F$ be globally generated with $h^0(\mathcal F)=2$ and $F\in |\mathcal F|$.
 Assume that
\begin{enumerate}
\item[(1)] $h^1(\mathcal M^2(F))=0$;
\item[(2)] $h^1(\mathcal M^2)+1\le h^1(\mathcal M^2(-F))$;
\item[(3)] $H^0(\mathcal M)\otimes H^0(F)\to H^0(\mathcal M(F))$ is surjective.
\end{enumerate}
 Then $\mathcal L(:=\mathcal M(F))$ fails to be normally generated, especially
 the 2-normality of $\phi_{\mathcal L}(C)$ does not hold.
\end{thm}

Note that hypotheses except $(3)$ in the above theorem could be
easily checked and possibly hold. For example, if $\mathcal M^2$ is
special then hypothesis $(2)$ naturally holds by $h^0 (\mathcal F) =
2, ~h^1(\mathcal M^2(F))=0$ and the base point free pencil trick.
For condition (3) we give a geometric criterion. The line bundle
$\mathcal L$ becomes extremal if $\mathcal L$ computes the Clifford
index of $C$. Generally it is hard to determine whether $\mathcal L$
computes the Clifford index of $C$ or not.

The following theorem makes it possible to find out line bundles
$\mathcal M$ and $\mathcal F$ satisfying the assumptions of Theorem
\ref{main2} such that $\mathcal M(F)$ is an extremal line bundle
with $h^1(\mathcal M(F))=2$.

\begin{thm}\label{thm3.3}
 Let $C$ be a smooth curve of genus $g$ and Clifford index
$c$ with $g\geq 2c+5$ which is neither hyperelliptic nor bielliptic.
If a line bundle $\mathcal M$ computes the Clifford index of $C$
with $(3c/2)+3 < {\deg} \mathcal M\leq g-1$, then $g=2c+5$ and
$|\mathcal M| =|{F +F'}|$ such that $|{F}|$ and $|F'|$ are penciles
of degree $c+2$. Moreover, $\mathcal M(F)$ computes the Clifford
index of $C$.
\end{thm}

We also show that the line bundle $\mathcal M$ is half canonical
unless $C$ is a $\frac{c+2}{2}$-fold covering of an elliptic curve
in Proposition \ref{cor3.2}. In fact, such a pair $(C,\mathcal M)$
can be found on a $K3$-surface whose Picard group is generated by a
hyperplane section in $\mathbb P^r$.

\begin{thm}\label{thm3.6}
Let $X$ be a general $K3$ surface in $\mathbb P^r$ $(r\ge 3)$ whose
Picard group is generated by a hyperplane section $H$ with ${\deg}
X=2r-2$. Let $C$ be a smooth irreducible curve on $X$ contained in
the linear system $|2H|$. Then
\begin{enumerate}
\item[(1)] $g=2c+5$ and there are pencils $|F|$, $|F'|$ of degree $c+2$
such that $|\mathcal O_C(1)|=|F+F'|$,

\item[(2)] both $\mathcal O_C(1)$ and $\mathcal O_C(1)(F)$ compute the
Clifford index of $C$,

\item[(3)] $\mathcal O_C(1)$ is normally generated, but $\mathcal
O_C(1)(F)$ is not.
\end{enumerate}
\end{thm}

Consequently, for any $r$ one can find couples $(C,\mathcal L)$ such
that $\mathcal L$ is an extremal line bundle on $C$ with
$h^1(\mathcal L)=2$ and $h^0(\mathcal L)=r+1$. For a curve $C$ of
genus $g\equiv 1(\mbox{mod} ~4)$, $C$ has both a normally generated
line bundle $\mathcal M$ and a non-normally generated line bundle
$\mathcal L$ with $h^1(\mathcal L)=2$ computing the Clifford index
of $C$ at the same time.

Recall the following: The Clifford index of the curve $C$ is defined
by
$$\cli(C):= \min \{ \cli(\mathcal L) : h^0 (\mathcal L) \geq 2, \ h^1 (\mathcal L) \geq 2\}.$$

A line bundle $\mathcal L$ is said to compute the Clifford index of
$C$ if $\cli(\mathcal L)=\cli(C)$ with $h^0(\mathcal L)\ge 2$ and
$h^1(\mathcal L)\ge 2$. The Clifford dimension of $C$ is defined by
$$r(C):= \min\{ h^0(\mathcal L)-1 : \mathcal L \text{ computes the
Clifford index of } C\}.$$ It has known that a general $k$-gonal
curve has Clifford dimension 1
(\cite{B}, \cite{Kek}). A smooth curve $C$ is called an {\it exceptional} curve if $r(C)\geq 2$.\\

\noindent{\bf Notations}

We denote the canonical line bundle on $C$ by $K$, $H^i(C,\mathcal
L)$ by $H^i(\mathcal L)$ and $H^i(C,\mathcal O (D))$ by $H^i (D)$.
We abuse the notations as follows: $L\in |\mathcal L|$, $|{\mathcal
O}(D)|=|D|$. For a divisor $D$ on $C$, we denote $\langle
D\rangle_{\mathcal L} $ the linear space spanned by $D$ in the
embedding associated to a very ample line bundle $\mathcal L$.

\section{Line bundles computing the Clifford index}

 In this section, we have the following results:
 Let $C$ be a curve of genus $g$ and Clfford index $c$
 which is neither hyperelliptic nor bielliptic, and a line bundle $\mathcal M$
 compute the Clifford index of $C$ with ${\deg} \mathcal M\leq g-1, ~h^0(\mathcal M)\ge 4$.
 Then there is a quadric hypersurface of rank $\leq 4$ containing
 $\varphi _{\mathcal M}(C)$
 if and only if $|\mathcal M| = |g^1_{c+2} + h^1_{c+2}|$ and $g=2c+5$.
 Thus if $g\ge 2c+5$ and ${\deg} \mathcal M> \frac{3c}{2}+3$, then $|\mathcal
 M|=|g^1_{c+2}+h^1_{c+2}|$ and $g=2c+5$. Moreover in
 this case, $\mathcal M$ is half-canonical unless
 $C$ is a $\frac{c+2}{2}$-fold covering of an elliptic curve.

 This is an extension of Corollary 3.2.5 in \cite{CM}.
 They showed that ${\deg} \mathcal M\le \frac{3c}{2}+3$ for
 a line bundle $\mathcal M$ computing the Clifford index of $C$ if $g>2c+4$ (resp.
 $g>2c+5$)
 and $c$ is odd (resp. even). In fact, the above is an extended result for the case
 $g=2c+5$ and $c$ is even. This is also comparable to the following:
 Let $C$ be an exceptional curve of genus $g$ and Clifford index
  $c$ and a line bundle $\mathcal M$ compute the Clifford dimension of
  $C$. Then $\varphi _{\mathcal M} (C)$ is
  not contained any quadric hypersurface of rank $4$ or less if $r(C)\geq
  3$. If $\deg \mathcal M> (3c/2)+3$, then $g=2c+4$ and $\mathcal M$ is
  half-canonical.
  Consult \cite{ELMS} for details .

\begin{Prop}[\cite{KKM}]\label{prop2.1}
Let $C$ be a  smooth  curve of genus $g $ and a line bundle
$\mathcal M$ compute the Clifford
 index of $C$ with ${\deg} \mathcal M\leq g-1$ and $h^0(\mathcal M)\geq 4$.
Then $\mathcal M$ is birationallly very ample unless $C$ is
hyperelliptic or bielliptic.
\end{Prop}

We note that the proposition also holds for such a linear system of
any degree, since the condition $d\le g-1$ is not used in the proof
of the proposition.

\begin{lem}\label{lem2.2}
Let $C$ be a smooth curve of genus $g$ and a line bundle $\mathcal
M$ compute the Clifford index of $C$ with ${\deg} \mathcal M \leq
g-1$. If $C$ has a base point free linear system $|F|$ with
$h^0(F)\geq 2$ and $h^0(\mathcal M(-F)) \geq  2$, then the linear
system $|F|$, $|\mathcal M(-F)|$ and $|\mathcal M(F)|$ compute the
Clifford index of $C$ .
\end{lem}
\begin{proof} We set ${\deg} F=f$. Since $\mathcal M$ computes the
Clifford index of $C$, we have $\cli (\mathcal M(-F))\geq
\cli(\mathcal M)$ and so $h^0(\mathcal M) \geq h^0 (\mathcal M(-F))
+ \frac{f}{2}\geq \frac{f}{2} +2$ for $h^0 (\mathcal M(-F))\geq 2$.
Assume $h^1 (\mathcal M(F))\leq 1$. Then by Rieman-Roch Theorem,
$h^0(\mathcal M(F))\leq f+1$ for ${\deg} \mathcal M\leq g-1$. Thus
$h^0(\mathcal M(F))- h^0(\mathcal M)\leq (f+1) -(\frac{f}{2} +2)\leq
\frac{f}{2}-1$. Then by the base point free pencil
trick(\cite{ACGH}, p 126), we get $\frac{f}{2}-1\geq h^0(\mathcal
M(F))- h^0(\mathcal M) \geq h^0(\mathcal M) - h^0(\mathcal
M(-F))\geq \frac{f}{2}$, which is a contradiction. Thus
$h^1(\mathcal M(F))\geq 2$ and so $\cli(\mathcal
M(F))\geq\cli(\mathcal M)$ which gives $h^0(\mathcal M(F))-
h^0(\mathcal M) \leq \frac{f}{2}$. Hence by the base point free
pencil trick, we have $\frac{f}{2} \geq h^0(\mathcal M(F))-
h^0(\mathcal M) \geq h^0(\mathcal M) - h^0(\mathcal M(-F))\geq
\frac{f}{2}$. Consequently, $\cli(\mathcal M(F))=\cli(\mathcal
M(-F))=\cli(\mathcal M)=\cli(C)$. Moreover all of them compute the
Clifford index of $C$, since $h^1(\mathcal M(F))\geq 2$ and
${\deg}(\mathcal M(-F))\leq g-1$. Set $|\mathcal M(-F)| =|F'|$. Then
$|F'|$ is base point free, since it computes the Clifford index of
$C$. Hence by the same argument for $|F'|$ instead of $|F|$, the
linear system $|F|=|\mathcal M(-F')|$ also computes the Clifford
index of $C$.
\end{proof}

For divisors $M$ and $E$ on $C$, let $(M,E)$ denote the greatest
common divisor of them.

\begin{lem} \label{lem2.3} Let a line bundle $\mathcal M$ compute the
Clifford index of $C$ with ${\deg} \mathcal M\le g-1$, $h^0(\mathcal
M) \ge 3$ and $\mathcal E$ be a line bundle with $h^0(\mathcal E)
\geq 2$. Then for any $P\in C$ there are divisors $M\in |\mathcal
M|$ and $E\in |\mathcal E|$ such that $(M,E)= P$ or $P+Q$ for some
$Q\in C$.
\end{lem}

\begin{proof} Let $P$ be an arbitrary point of $C$ and  $E$ a
divisor in $|\mathcal E|$ containing $P$. Set $E= P+\Sigma P_i $.
Let $B$ be the base locus of $|\mathcal M(-P)|$. Then $B$ is either
zero divisor or degree one divisor $Q$ for some $Q\in C$, since
$\mathcal M$ computes the Clifford index of $C$ with $h^0(\mathcal
M) \geq 3$. We set $\mathcal G: =\mathcal M(-P -B)$. Then $\mathcal
G$ is base point free of $h^0(\mathcal G)\ge 2$. Hence there is a
$G\in |\mathcal G|$ such that $(G, \Sigma P_i)=$ zero divisor. Thus
$(M,E)= P$ or $P+Q$ for $M=G +P+B \in |\mathcal M|$.
\end{proof}

Using the above results, we get the following proposition.

\begin{Prop} \label{prop3.1}
 Let $C$ be a smooth curve of genus $g$ and Clifford index
$c$ which is neither hyperelliptic nor bielliptic and let $\mathcal
M$ be a line bundle
  computing the Clifford index of $C$ with ${\deg} \mathcal M\leq g-1$
  and $h^0(\mathcal M) \geq 4$. Then there is a quadric hypersurface of rank $\leq
  4$ containing $\varphi _{\mathcal M} (C)$ if and only if $|\mathcal M| =|F +F'|$ and
  $g=2c+5$, where $|F|$ and $|F'|$ are pencils of degree $c+2$.
  In this case, $\mathcal M(F)$ computes the Clifford index of $C$.
\end{Prop}

\begin{proof} By the assumption and
Proposition \ref{prop2.1}, the morphism $\varphi _{\mathcal M}$ is
birational. Fix a quadric hypersurface $Q$ of rank $\leq 4$
containing $\varphi _{\mathcal M} (C)$. Let $|F_1|$ and $|F_2|$ be
the complete linear systems induced by two pencils on $Q$. Then
$|\mathcal M| =|F_1 +F_2|$. Let $|F|$ be a base point free pencil
which is a subsystem  of $|F_1| $. Then $h^0 (\mathcal M(-F))\geq
h^0 (F_2) \geq 2$. By Lemma \ref{lem2.2}, $|F| =g^1_{c+2}$, both
$|\mathcal M(-F)|$ and $|\mathcal M (F)|$ also compute the Clifford
index of $C$. In particular, we have $h^0(K\otimes\mathcal
M(F)^{-1})\geq 2.$

\noindent{Claim}: $|\mathcal M(-F)|=g^1_{c+2}$

To prove this, we assume $h^0(\mathcal M(-F))\geq 3$. Let $P$ be an
arbitrary point of $C$. By Lemma \ref{lem2.3}, there are divisors
$G_1 \in |\mathcal M(-F)|$ and $G_{2} \in |K\otimes\mathcal
M(F)^{-1}|$ such that $(G_1,G_2)=P$ or $P+R$ for some $R\in C$. Thus
if we set $M=G_1 +F$ and $E=G_2+F$, then $M\in |\mathcal M|$ and
$E\in |K\otimes\mathcal M^{-1}|$ with $(M,E)=F+P$ or $F+P+R$. First
we assume $(M,E)=F+P$. Then by the sheaf exact sequence (see
\cite{H}, p 345)
$$0\rightarrow {\mathcal O}((M,E)) \rightarrow {\mathcal O}(M) \oplus {\mathcal O} (E)
\rightarrow {\mathcal O} (K(-(M,E)))\rightarrow 0,$$ we have
$\cli(F+P)\leq \cli(M)$. It is impossible because
$\cli(M)=\cli(F)=c$.

As a consequence, $(M,E)=F+P+R$ for some $R\in C$. Using the above
sheaf exact sequence, $\cli(F+P+R)=\cli(M)=\cli(F)$, and hence
$h^0(K(-F-P-R))=h^0(K(-F))-1$. On the other hand, we have
$h^0(K(-F))\ge h^0(\mathcal M)\ge 4$. Thus by Proposition
\ref{prop2.1}, $|K(-F)|$ is birationally very ample since
$\cli(K(-F))=c$. Thus $h^0(K(-F-P-R))=h^0(K(-F))-2$ except only
finite pairs $(P,R)$. It contradicts to the arbitrary choice of the
point $P$. Thus $h^0(\mathcal M(-F)) =2$ and so $|{\mathcal
M(-F)}|=g^1_{c+2}$.

If we exchange the roles of $|{\mathcal M(-F)}|$ and
$|K\otimes\mathcal M(F)^{-1}|$ in the claim, then we have
$h^0(K\otimes\mathcal M(F)^{-1})=2$ and so $|K\otimes\mathcal
M(F)^{-1}|=g^1_{c+2}$ since ${K\otimes\mathcal M^{-1}}$ also
computes the Clifford index of $C$. Thus both $|\mathcal M|$ and
$|K\otimes\mathcal M^{-1}|$ are sums of two linear pencils of degree
$c+2$, which proves the result.
\end{proof}

\noindent {\bf Proof of Theorem \ref{thm3.3}.} The condition
$(3c/2)+3 < {\deg} \mathcal M$ yields that $h^0(\mathcal M)\ge 3$.
If $h^0(\mathcal M) =3$, then we have $c=1$ since $\mathcal M$
computes the Clifford index of $C$ and $(3c/2)+3 < {\deg} \mathcal M
=c+4$. Thus $C$ is a plane quintic, which cannot occur since $g\geq
2c+5$. Accordingly, $h^0(\mathcal M) \geq 4$.

Denote $h^0(\mathcal M) =r+1$. Suppose that the result does not
hold, then by Proposition \ref{prop3.1}, the image curve $C'
=\varphi _{\mathcal M} (C)$ is not contained any quadric
hypersurface of rank$\leq 4$. Then by the sheaf exact sequence
$$0\rightarrow {\mathcal I}_{C'}(2) \rightarrow {\mathcal O}_{{\mathbb P}^r}(2)\rightarrow {\mathcal O}_{C'}(2)
\rightarrow 0,$$ we have  $h^0 (\mathcal M^2)\geq 4r-2$. Hence
$\cli(\mathcal M^2)\leq 2c-4r+6$ and $h^1 (\mathcal M^2)\geq 2$
since ${\deg} \mathcal M= c+2r$ and $g\geq 2c+5$. Accordingly,
$\cli(\mathcal M^2)\geq c$ and so $4r\leq c+6$. Thus we have ${\deg}
\mathcal M=c+2r\le
c+\frac{c}{2}+3$ which is a contradiction. \qed\\

We can show that such a line bundle $\mathcal M$ is generally
half-canonical for the boundary case $g=2c+5$. It is comparable to
the line bundle computing the Clifford dimension on an exceptional
curve $C$ of genus $g(C)=2\cli(C)+4$.

\begin{Prop}\label{cor3.2} Let $C$ and $\mathcal M$ be the same as Theorem
\ref{thm3.3}. Assume there is a quadric hypersurface of rank $\leq
4$ containing $\varphi _{\mathcal M} (C)$. Then $\mathcal M$ is
half-canonical unless $C$ is a $\frac{c+2}{2}$-fold covering of an
elliptic curve.
\end{Prop}

\begin{proof}
By the above theorem we may set $|\mathcal M|=|F+F_1|$ and
$|K\otimes\mathcal M^{-1}|=|F+F_2|$ such that $|F|,|F_1|$ and
$|F_2|$ are base point free pencils of degree $c+2$. Assume that
$\mathcal M$ is not half-canonical, i.e., $|F_1|\neq |F_2|$.
Consider $\phi_{F_1}\times \phi_{F_2}: C\rightarrow \mathbb P^3$ and
let $C'$ be a smooth model of $\phi_{F_1}\times\phi_{F_2}(C)$. Then
we have a morphism $\psi:C\to C'$ such that
$\pi\circ\psi=\phi_{F_1}\times\phi_{F_2}$ where $\pi$ is a
normalization morphism from $C'$ to $\phi_{F_1}\times\phi_{F_2}(C)$.
Take a divisor $G_i$ on $C'$ such that $F_i=\psi^*(G_i)$ and let
$m:={\deg} \psi$, then $|G_i|$'s are base point free pencils of
degree $\frac{c+2}{m}$ on $C'$ since $|F_i|$'s are base point free.
Then we have the following commutative diagram:

\begin{picture}(300,100)

\put(80,5){$C$}

\put(100,5){\vector(1,0){80}}

\put(95,20){\vector(3,2){100}}

\put(205,90){$C'$}

\put(210,80){\vector(0,-1){60}}

\put(215,53){$\pi=\phi_{G_1}\times \phi_{G_2}$}

\put(132,53){$\psi$}

\put(120,13){$\phi_{F_1}\times \phi_{F_2}$}

\put(186,5){$\phi_{F_1}\times \phi_{F_2}(C)\subset \mathbb P^3$}

\end{picture}

If we note that $\phi_{G_1}\times \phi_{G_2}$ is birational then
there exists $G_i\in |G_i|$ such that $(G_1,G_2)=Q$ for any $Q\in
C'$. Fix $G_1:=Q_1+\cdots+Q_{\frac{c+2}{m}}$ and let
$R_i:=\psi^*(Q_i)$, then
$\psi^*(G_1)=R_1+\cdots+R_{\frac{c+2}{m}}\in |F_1|$. Hence
$F+\psi^*(G_1)\in |\mathcal M|$ and there exists $G_{2,i}\in |G_2|$
such that $(G_1,G_{2,i})=Q_i$ for each $i$. Let $M:=F+\psi^*(G_1)$
and $E_i:=F+\psi^*(G_{2,i})$. Then $(M,E_i)=F+\psi^*(Q_i)=F+R_i$.
Note that $\cli(F+R_i)=c$ for any $i$ by the short exact sequence
$$0\rightarrow {\mathcal O}(M,E_i) \rightarrow {\mathcal O}(M)
\oplus {\mathcal O} (E_i) \rightarrow {\mathcal O}
(K(-(M,E_i)))\rightarrow 0.$$ Let $s:=\frac{c+2}{m}$.

\noindent {\bf claim :} $\cli(F+R_1+\dots +R_k)=c$ for
$k=1,\dots,s$.

First we prove that $\cli(F+R_1+R_2)=c$. By the definition of
Clifford index, we obtain $\cli(F+R_1+R_2)={\deg}
(F+R_1+R_2)-2r(F+R_1+R_2)\ge c$.  So, we get ${\dim} \langle
F+R_1+R_2\rangle_K \ge c+m$ from the geometric Riemann-Roch Theorem.

On the other hand, ${\dim} \langle F+R_i\rangle_K = c+\frac{m}{2}$
and ${\dim} \langle F\rangle_K \cap \langle R_i \rangle_K =
\frac{m}{2}-1$ since $\cli(F+R_i)=c, ~{\dim}\langle F\rangle_K=c$
and ${\dim}\langle R_i\rangle_K=m-1$. It produces ${\dim} \langle
F+R_1+R_2\rangle_K \le c+m$. Therefore, $\cli(F+R_1+R_2)=c$. In the
same manner, one can prove that ${\dim} \langle
F+R_1+\dots+R_k\rangle_K = c+\frac{km}{2}$, which gives a proof of
the claim.

If $s>2$, then we lead to ${\deg}(F+R_1+\dots +R_{s-1}) >
\frac{3c}{2}+3$, which gives a contradiction to Theorem
\ref{thm3.3}. This contradiction gives us that $\mathcal M$ is
half-canonical. If $s=2$, then there are different pencils $|G_1|$,
$|G_2|$ of degree $2$ on $C'$. Thus $C'$ is an elliptic curve and
hence $C$ is a $\frac{c+2}{2}$-fold covering of the elliptic curve
$C'$. We have excluded this case. In all, we obtain the result.
\end{proof}

{\it\bf {Proof of Theorem \ref{thm3.6}}.} Since $X$ is a K3-surface
in $\mathbb P^r$, $\mathcal O_C(2)$ is the canonical bundle of $C$,
$g(C)=2c+5$ and ${\deg} \mathcal O_C(1)=2c+4$ from adjunction
formula. By Green and Lazarsfeld Theorem (\cite{GL1}, Theorem),
$\mathcal O_C(1)$ computes the Clifford index of $C$. Therefore the
results (1) and (2) follow from Theorem \ref{thm3.3}.

$X$ is projectively normal since a hyperplane section is a canonical
curve. From the exact sequence: $$0\to H^0(\mathcal I_X(2))\to
H^0(\mathcal O_{\mathbb P^r}(2))\to H^0(\mathcal O_X(2))\to 0,$$ one
can see that $h^0(\mathcal I_X(2))=\binom{r+2}{2}-{\deg} C-2$ by the
Riemann-Roch Theorem. We have $h^0(\mathcal
I_C(2))=\binom{r+2}{2}-{\deg} C-1$ by the short exact sequence:
$$0\to H^0(\mathcal I_X(2))\to H^0(\mathcal I_C(2))\to H^0(\mathcal
O_X)\to 0.$$  The morphism $\mu_m: H^0(\mathcal O_{\mathbb
P^r}(m))\to H^0(\mathcal O_C(m))$ is surjective for $m=2$ by
$h^0(\mathcal O_C(2))=g={\deg} C+1$. For $m\ge 3$, the surjectivity
of $\mu_m$ follows from the proof of Theorem 3.6 in \cite {ELMS}.
Therefore $\mathcal O_C(1)$ is normally generated.

Corollary \ref{cor4.4} in the next section implies that $\mathcal
O_C(1)(F)$ fails to be normally generated. \qed

\section{Extremal line bundles $\mathcal L$ with $h^1(\mathcal L)=2$}
We start this section to give a geometric interpretation of the
surjectivity of Theorem \ref{main2}.
\begin{lem}\label{empty}
 Let $\mathcal M$ be a very ample line bundle on a smooth curve $C$,
 $\mathcal F$ a globally generated line bundle with $h^0(\mathcal F)=2$ and $\mathcal L:=\mathcal M\otimes\mathcal F$.
 For any two distinct divisors $F, F'\in \mathcal F$,
 $\langle F\rangle_{\mathcal L}\cap \langle F'\rangle_{\mathcal L}=\emptyset$ if and only if the cup product morphism
 $\mu: H^0(\mathcal M)\otimes H^0(F)\to H^0(\mathcal L)$ is surjective.
\end{lem}
\begin{proof}
  By the base point free pencil trick, we have the following exact sequence;
 \begin{equation*}\label{bpf}
   0\to H^0(\mathcal M(-F))\to H^0(\mathcal M)\otimes H^0(F)\stackrel{\mu}\rightarrow H^0(\mathcal L).
 \end{equation*}
  Let $r:=h^0(\mathcal M)-1$. We have
  ${\dim} {\langle F\rangle_{\mathcal L}}=h^0(\mathcal L)- r -2$ and
  ${\dim} \langle F+F'\rangle_{\mathcal L}=h^0(\mathcal L)-h^0(\mathcal M(-F))-1$.

  If $\mu$ is surjective, then
  $h^0(\mathcal L)=2(r+1)-h^0(\mathcal M(-F))$. Therefore
  ${\dim} \langle F+F'\rangle_{\mathcal L}= 2{\dim} \langle F\rangle_{\mathcal L}+1$ which implies
   $\langle F\rangle_{\mathcal L}\cap \langle F'\rangle _{\mathcal L'}=\emptyset$.

  If $\langle F\rangle _{\mathcal L}\cap \langle F'\rangle_{\mathcal L}=\emptyset$, then
  ${\dim} \langle F+F'\rangle_{\mathcal L}= 2{\dim} \langle F\rangle_{\mathcal L}+1$.
  This yields $h^0(\mathcal L)=2(r+1)-h^0(\mathcal M(-F))$,
  so $\mu$ is surjective.
 \end{proof}

{\it\bf {Proof of Theorem \ref{main2}}.}
 First, we claim that $F$ fails to impose independent conditions
  on quadrics in $\langle F\rangle_{\mathcal M}$, i.e. $H^1(\mathcal I_{F/\langle F\rangle_{\mathcal M}}(2))\neq 0$.
  Since  $H^2(\mathcal I_{F/\mathbb P^m}(2))=0$,
  we have the following exact sequence:
  \begin{equation*}
   H^1(\mathcal I_{F/\mathbb P^m}(2))\to
   H^1(\mathcal I_{F/C}(2))\to
   H^2(\mathcal I_{C/\mathbb P^m}(2)) \to 0
  \end{equation*}
  where $\mathbb P^m:=\mathbb P(H^0(\mathcal M))$.
  Note that we have $H^2(\mathcal I_C(2))\simeq H^1(\mathcal O_C(2))\simeq H^1(\mathcal M^2)$ and
  $H^1(\mathcal I_{F/C}(2))\simeq H^1(\mathcal M^2(-F))$. By the assumption (2), we conclude
  that $H^1(\mathcal I_{F/\mathbb P^m}(2))\neq 0$. Since $H^1(\mathcal I_{F/\mathbb P^m}(2)) =
  H^1(\mathcal I_{F/\langle F\rangle_{\mathcal M}}(2))$, $F$ fails to impose independent conditions
  on quadrics in $\langle F\rangle_{\mathcal M}$.

 By Lemma \ref{empty} hypothesis (3) makes it possible to take
 two divisors $F,F'\in |g^1_{c+2}|$ such that
 $\langle F\rangle_{\mathcal L}\cap \langle F'\rangle_{\mathcal L}=\emptyset$ for $\mathcal L:=\mathcal M(F)$.
 Consider the projection $\pi_{F'}$ of $\phi_{\mathcal L}(C)$ from $F'$.
 Then we have the following commutative diagram:

\begin{picture}(300,100)

\put(80,75){$C$}

\put(95,80){\vector(1,0){85}}

\put(90,65){\vector(2,-1){90}}

\put(186,75){$\phi_{\mathcal L}(C)\subset \mathbb P^r(:=\mathbb
P(H^0(\mathcal L)))$}

\put(190,70){\vector(0,-1){40}}

\put(200,45){$\pi_{F'}$: projection}

\put(130,87){$\phi_{\mathcal L}$}

\put(115,33){$\phi_{\mathcal M}$}

\put(186,15){$\phi_{\mathcal M}(C)\subset \mathbb P^m(:=\mathbb
P(H^0(\mathcal M)))$}

\end{picture}


 Since $\langle F\rangle_{\mathcal L}\cap \langle F'\rangle_{\mathcal
 L}=\emptyset$, we see that $H^1(\mathcal I_{F/\langle F\rangle_{\mathcal M}}(2))\neq
 0$ if and only if $H^1(\mathcal I_{F/\langle F\rangle_{\mathcal L}}(2))\neq
 0$. Therefore $F$ fails to impose independent conditions on
 quadrics in $\langle F\rangle_{\mathcal L}$, i.e. $H^1(\mathcal I_{F/\langle F\rangle_{\mathcal L}}(2))\neq 0$.
 This is equivalent to $H^1(\mathcal I_{F/\mathbb P^r}(2))\neq 0$.
 From the following exact sequence
 \begin{equation*}
 0\to\mathcal I_{C/\mathbb P^r}(2)\to\mathcal I_{F/\mathbb P^r}(2)\to
 \mathcal I_{F/C}(2)\to 0,
 \end{equation*}
 one can see that $H^1(\mathcal I_{C/\mathbb P^r}(2))\neq 0$
 since  $H^1(\mathcal I_{F/C}(2))=H^1(\mathcal L^2(-F))=H^1(\mathcal M^2(F))=0$
 by the hypothesis $(1)$.
 Thus $\mathcal L(:=\mathcal M(F))$ fails to be normally generated.\qed

  \begin{rmk}
    If ${\deg} \mathcal L\ge g+1$, then the multiplication map
    $$
    \mu_m:H^0(\mathcal L)\otimes H^0(\mathcal L^m)\to H^0(\mathcal L^{m+1})
    $$
    are surjective for $m\ge 2$ $($\cite{Gr},~ Theorem $(4.e.1))$.
    Therefore, $\mathcal L$ is normally generated if $\mathcal L$
    satisfies $2$-normality.
  \end{rmk}

\begin{cor}\label{cor4.4} Let $C$, $\mathcal M$ and $F$ be the same
as in Theorem \ref{thm3.3}. Then $\mathcal M(F)$ is an extremal line
bundle with $h^1(\mathcal M(F))=2$.
\end{cor}

\begin{proof}
The hypothesis $(1)$ of Theorem \ref{main2} holds trivially since
${\deg} \mathcal M=g-1$. From Lemma \ref{lem2.2}, one can see that
$K\otimes\mathcal M^{-2}(F)=K\otimes\mathcal M^{-1}(-F')$. Thus
$h^1(\mathcal M^2(-F))=h^1(\mathcal M(F'))=2$. Also we have
$h^1(\mathcal M^2)\le 1$ since ${\deg} \mathcal M=g-1$. Therefore
the hypothesis $(2)$ of the Theorem \ref{main2} is satisfied. Using
the base point free pencil trick and the Riemann-Roch Theorem
$\mu:H^0(\mathcal M)\otimes H^0(F)\rightarrow H^0(\mathcal M(F))$ is
surjective, so the hypothesis (3) holds. Hence the result follows
from Theorem \ref{main2} and Theorem \ref{thm3.3}.
\end{proof}

\section{Examples and questions}

In this section, we observe examples of extremal line bundles
$\mathcal L$ with $h^1(\mathcal L)\le 2$.

\begin{Prop}[\cite{GL}, Remark 2.6]
Let $C$ be a $k$-gonal curve such that $\cli(C)$ is computed by a
pencil $g^1_k$. Let $\sum^4_{i=1}P_i$ be a divisor on $C$ such that
$h^0(g^1_k-P_i-P_j))=0$ for all $i,j\in \{1,\cdots,4\}$. Then
$\mathcal L:=K(-g^1_k+\sum^4_{i=1} P_i)$ is a nonspecial extremal
line bundle on $C$.
\end{Prop}

The following examples also can be found in \cite{Ko}.

  \begin{Prop}
     Let $C$ be a nonsingular plane curve degree $d\ge 5$ and $H$ a
     line section of $C$. Take a degree 4 divisor $Z\le H$
     and put $D=H-Z$.
     Then $K(-D)$ is an extremal line bundle on $C$ with $h^1(K(-D))=1$.
  \end{Prop}

  \begin{proof}
     Since $C$ has no pencil of degree $\le d-2$, $K(-D)$ is very ample.
     By the Riemann-Roch theorem $Z$ spans a line in $\mathbb P^r:=\mathbb P(H^0(K(-D)))$,
     and hence fails to impose independent conditions on quadrics,
     i.e. $h^1(\mathcal I_{Z/\mathbb P^r}(2))\neq 0$.
     Now consider the following exact sequence;
     \begin{equation*}
      0\to\mathcal I_{C/\mathbb P^r}(2)\to\mathcal I_{Z/\mathbb P^r}(2)\to
      \mathcal I_{Z/C}(2)\to 0
     \end{equation*}
     Since  $H^1(\mathcal I_{Z/C}(2))=H^1((K(-D))^2(-Z))=0$,
     $H^1(\mathcal I_{C/\mathbb P^r}(2))\neq 0$, i.e.
     $K(-D)$ fails to be normally generated.
     Since $\cli(K(-D))=d-4$, $K(-D)$ is an extremal line bundle.
  \end{proof}

  Let $C$ be an exceptional curve with $g=2c+4$ and $\mathcal M$ compute
  the Clifford dimension of $C$. Then we have ${\deg} \mathcal
  M=4r-3$, $\cli(C)=2r-3$ and $g=4r-2$ where $h^0(\mathcal M)=r+1$. It is
  known that $\phi_{\mathcal M}(C)$ is projectively normal and
  has a $(2r-3)$-secant $(r-2)$-space divisor
  $D$ (\cite{CM}, Theorem A). Konno proved that the
  line bundle $K(-D)$ is extremal with $h^1(K(-D))=1$ in \cite{Ko}.
 We reprove it using Theorem \ref{main2}.

  \begin{Prop}\label{cor4.3}
    Let $(C,\mathcal M, D)$ be as above. Then
    $K(-D)$ is an extremal line bundle with $h^1(K(-D))=1$.
  \end{Prop}

  \begin{proof}
   Let $\mathcal F\cong \mathcal M(-D)$, then $K(-D)\cong \mathcal M\otimes \mathcal
   F$ since $\mathcal M$ is a half-canonical.
  Conditions (1), (2) of Theorem \ref{main2} clearly hold.
   By the base point free pencil trick we have the following:
   $$0\to H^0(\mathcal M\otimes\mathcal F^{-1})\to H^0(\mathcal M)\otimes
   H^0(\mathcal F)\stackrel{\mu}\rightarrow H^0(\mathcal M\otimes \mathcal F).$$
   Then $\mu$ is surjective since $h^0(\mathcal M\otimes \mathcal F)=h^0(K(-D))=c+4$,
   $h^0(\mathcal M\otimes\mathcal
   F^{-1})=1$ and $h^0(\mathcal M)=\frac{c+5}{2}$. Theorem \ref{main2} implies
   that $K(-D)$ is an extremal line bundle
   since $\cli(K(-D))=\cli(C)$.
  \end{proof}

Even though the curves in next examples are lying on K3 surfaces, we
do the works with explicit calculations.

 \begin{ex}
  Let $C$ be a smooth complete intersection of smooth surfaces of degree 2 and
  4 in $\mathbb P^3$. Then ${\deg} C=8$, $g(C)=9$, and $gon(C)=4$.
  Therefore H. Martens' Theorem implies $\cli(C)=2$. For
  $\mathcal M:=\mathcal O_C(1)$, we have $|\mathcal M|=|F+F'|$
  where $F$ is a pencil of degree 4 by
  Theorem \ref{thm3.3}.
  Thus $\mathcal M(F)$ fails to be normally generated by Theorem \ref{main2}. So $\mathcal M(F)$ is
  extremal with $h^1(\mathcal M(F))=2$.
 \end{ex}

\begin{ex}\label{ex4.7}
 Let $C$ be a smooth complete intersection of smooth hypersurfaces of degree $2,2,3$
 in $P^4$. Then there is a quadric hypersurface $Q$ of rank $\le 4$ and $\mathcal M(F)$
 is an extremal line bundle with $h^1(\mathcal M(F))=2$. Here $\mathcal M=\mathcal O_C(1)$
 and $F$ is given by a ruling of $Q$ of degree 6.
\begin{proof}
 Note that $g(C)=13$. By Lazarsfeld's theorem $($\cite{La}, Example 4.12$)$,
 $C$ can not have a pencil $g^1_5$. It implies
 that Cliff($C)\ge 3$ by Coppens and Martens' Theorem $($\cite{CM}$)$.
 Assume that Cliff($C)=3$. If the Clifford dimension $r(C)$ of $C$ is 2, then $C$ is isomorphic
 to a smooth plane septic, which cannot occur since $g(C)=13$. Therefore $r(C)\ge
 3$ and there is a $g^{r(C)}_d$ with
 $d=3+2r(C)\le 12=g(C)-1$, and so $r(C)\le 4$. Hence $C$ is an ELMS
 curve (\cite{ELMS} Section 5), which is a contradiction because of $g(C)\neq 2c+4$.
 Thus $\cli(C)=4$ since $|\mathcal M|=g^4_{12}$ where $\mathcal
 M:=\mathcal O_C(1)$.
 Note that $C$ is always contained in a quadric of rank $4$ or less because
 $C$ is contained in two quadrics in $\mathbb P^4$.
Whence Proposition \ref{prop3.1} yields $|\mathcal M|=|F+F'|$ such
that $|F|, |F'|$ are base point free pencils of degree $c+2$. By
Theorem \ref{main2}, $\mathcal M(F)$ fails to be normally generated.
So $\mathcal M(F)$ is extremal with $h^1(\mathcal M(F))=2$.
 \end{proof}
\end{ex}

\begin{ex}
 Let $C\subset \mathbb P^5$ be a smooth complete intersection
 of four quadric hypersurfaces.
 Then there is a quadric hypersurface $Q$ of rank $\le 4$ and $\mathcal M(F)$
 is an extremal line bundle with $h^1(\mathcal M(F))=2$. Here $\mathcal M=\mathcal O_C(1)$
 and $F$ is given by a ruling of $Q$ of degree 8.
\end{ex}

\begin{proof}
 As in example \ref{ex4.7}, one can show that $g(C)=17$ and
 Cliff$(C)=6$. Note that the quadrics of rank $\le 4$
 in $\mathbb P^5$ form a closed subvariety of codimension $3$
 in the projective space of all quadrics in $\mathbb P^5$. Also,
 $C$ is always contained in a quadric of rank $\le 4$ because
 $C$ is contained in four quadrics.
 Therefore, by Proposition \ref{prop3.1},
 $|\mathcal M|=|F+F'|$ for any hyperplane section $\mathcal M$.
 By Theorem \ref{main2}, $\mathcal M(F)$ fails to be normally
 generated. So $\mathcal M(F)$ is extremal with $h^1(\mathcal M(F))=2$.
\end{proof}

Finally we ask a couple of questions relating to the above results.

\begin{question}
Is there an extremal line bundle $\mathcal L$ on a smooth curve with
$h^1(\mathcal L)\ge 3$?
\end{question}

\begin{question}
Can we find a smooth curve which does not lie on a $K3$ surface, but
has an extremal line bundle $\mathcal L$ with $h^1(\mathcal L)\ge
2$?
\end{question}

\end{document}